\documentclass[12pt,reqno]{amsart}
\usepackage{amssymb, latexsym}
\usepackage{url}
\usepackage{a4wide}

\begin{document}
 \bibliographystyle{plain}

 \newtheorem{theorem}{Theorem}
 \newtheorem{lemma}{Lemma}
 \newtheorem{corollary}{Corollary}
 \newtheorem{problem}{Problem}
 \newtheorem{conjecture}{Conjecture}
 \newtheorem{definition}{Definition}
 \newcommand{\mc}{\mathcal}
 \newcommand{\rar}{\rightarrow}
 \newcommand{\Rar}{\Rightarrow}
 \newcommand{\lar}{\leftarrow}
 \newcommand{\lrar}{\leftrightarrow}
 \newcommand{\Lrar}{\Leftrightarrow}
 \newcommand{\zpz}{\mathbb{Z}/p\mathbb{Z}}
 \newcommand{\mbb}{\mathbb}
 \newcommand{\A}{\mc{A}}
 \newcommand{\B}{\mc{B}}
 \newcommand{\cc}{\mc{C}}
 \newcommand{\D}{\mc{D}}
 \newcommand{\E}{\mc{E}}
 \newcommand{\F}{\mc{F}}
 \newcommand{\G}{\mc{G}}
  \newcommand{\ZG}{\Z (G)}
 \newcommand{\FN}{\F_n}
 \newcommand{\I}{\mc{I}}
 \newcommand{\J}{\mc{J}}
 \newcommand{\M}{\mc{M}}
 \newcommand{\nn}{\mc{N}}
 \newcommand{\qq}{\mc{Q}}
 \newcommand{\U}{\mc{U}}
 \newcommand{\X}{\mc{X}}
 \newcommand{\Y}{\mc{Y}}
 \newcommand{\itQ}{\mc{Q}}
 \newcommand{\C}{\mathbb{C}}
 \newcommand{\R}{\mathbb{R}}
 \newcommand{\N}{\mathbb{N}}
 \newcommand{\Q}{\mathbb{Q}}
 \newcommand{\Z}{\mathbb{Z}}
 \newcommand{\ff}{\mathfrak F}
 \newcommand{\fb}{f_{\beta}}
 \newcommand{\fg}{f_{\gamma}}
 \newcommand{\gb}{g_{\beta}}
 \newcommand{\vphi}{\varphi}
 \newcommand{\whXq}{\widehat{X}_q(0)}
 \newcommand{\Xnn}{g_{n,N}}
 \newcommand{\lf}{\left\lfloor}
 \newcommand{\rf}{\right\rfloor}
 \newcommand{\lQx}{L_Q(x)}
 \newcommand{\lQQ}{\frac{\lQx}{Q}}
 \newcommand{\rQx}{R_Q(x)}
 \newcommand{\rQQ}{\frac{\rQx}{Q}}
 \newcommand{\elQ}{\ell_Q(\alpha )}
 \newcommand{\oa}{\overline{a}}
 \newcommand{\oI}{\overline{I}}
 \newcommand{\dx}{\text{\rm d}x}
 \newcommand{\dy}{\text{\rm d}y}
\newcommand{\cal}[1]{\mathcal{#1}}
\newcommand{\cH}{{\cal H}}
\newcommand{\diam}{\operatorname{diam}}

\parskip=0.5ex

\title[]{The Duffin-Schaeffer Conjecture\\ with  extra divergence }

\author{Alan Haynes, Andrew Pollington and  Sanju Velani}

\begin{abstract}
Given a nonnegative function $\psi : \N \to \R $, let  $W(\psi)$
denote the set of real numbers $x$ such that  $|nx -a| < \psi(n) $
for infinitely many reduced rationals $a/n \ (n>0) $. A consequence
of our main result is  that $W(\psi)$ is of full Lebesgue measure if
there exists an $\epsilon
> 0 $ such that
$$
\textstyle
\sum_{n\in\N}\left(\frac{\psi(n)}{n}\right)^{1+\epsilon}\!\!\!\varphi
(n)=\infty \, .
$$
The Duffin-Schaeffer Conjecture is the corresponding statement with
$\epsilon = 0$ and represents a fundamental unsolved problem in
metric number theory. Another consequence is that $W(\psi)$ is of
full Hausdorff dimension if the above sum with $\epsilon = 0$
diverges; i.e. the dimension analogue of the Duffin-Schaeffer
Conjecture is true.
\\[2ex]

\end{abstract}

\thanks{\!AH:~Research supported by EPSRC grant
EP/F027028/1.}

\thanks{AP:~Research supported by the NSF}

\thanks{SV:~Research supported by EPSRC grants EP/E061613/1 and EP/F027028/1.}

\allowdisplaybreaks \maketitle

\section{Introduction}
Throughout this note we will use the following standard notation
from elementary number theory: $p$ denotes a prime number, $\mu
(n)$ is the M\"{o}bius function, $\varphi (n)$ is the Euler phi
function, $\omega (n)$ denotes the number of distinct prime
divisors of $n$, and $\tau (n)$ is the number of positive integers
which divide $n$. Also we use $\lambda$ to denote Lebesgue measure
on $\R/\Z$ and $\dim X $ to denote the Hausdorff dimension of a
set $X$.

Let $\psi:\N\rar\R$ be a nonnegative arithmetical function and for
each positive integer $n$ define $\E_n\subseteq \R/\Z$ by
\begin{equation*}
\E_n \, := \,
\bigcup_{\substack{a=1\\(a,n)=1}}^n\left(\frac{a-\psi
(n)}{n},\frac{a+\psi (n)}{n}\right).
\end{equation*}
Denote the collection of points $x\in\R/\Z$ which fall in
infinitely many of the  sets $\E_n$ by $W(\psi)$. In other words,
$$
\textstyle{W(\psi):=
\limsup_{n\to\infty}\E_n:=\bigcap_{m=1}^\infty\ \bigcup_{n\ge
m}\E_n \ }  \ .
$$

\noindent The famous Duffin-Schaeffer Conjecture
\cite{DuffinSchaeffer} dates back to 1942 and is the assertion
that $\lambda (W(\psi))=1$ if
\begin{equation}\label{divcond1}
\sum_{n\in\N}\frac{\psi (n) }{n}  \, \varphi (n) =\infty \, .
\end{equation}
Equivalently, the $\limsup$ set $W(\psi)$  is of  full Lebesgue
measure if the sum of the measures of the sets  $\E_n$ diverges.
Although various partial results\footnote{A partial result of
particular importance is Khintchine's theorem 
from 1924. This fundamental theorem implies the Duffin-Schaeffer
Conjecture under the assumption that $\psi$ is monotonic.  This is
hardly surprising since the conjecture is a direct result of
attempting to remove monotonicity from Khintchine's theorem.} have
been established (see \cite{HarmanMNT} for details and references)
the full conjecture represents one of the most difficult and
profound unsolved problems in metric number theory. Our goal here is
to prove a  weaker statement in which `extra divergence' is assumed.
To this end,  define the function $f:[0,\infty)\rar\R$ by
\begin{equation*}
f(x) \, := \, \begin{cases}0&\text{ if }x=0,\\
x\cdot\exp\left\{\frac{\log x}{2\log(-\log x)}\right\}&\text{ if }0<x<1,\\
1&\text{ if }x\ge 1.\end{cases}
\end{equation*}
It is easily  seen that as $x$ tends to zero from above,  $f(x)$
tends to zero  faster than $x \,  (- \log x)^{-\epsilon}$ but more
slowly than $x^{1+\epsilon}$ for any $\epsilon>0$.

\begin{theorem}\label{extradivthm}
Let $\psi$ be any nonnegative arithmetical function and let $f$ be
defined as above. Then $\lambda (W(\psi))=1$ if
\begin{equation}\label{divcond2}
\sum_{n\in\N}f\left(\frac{\psi(n)}{n}\right)\varphi (n)=\infty \, .
\end{equation}
\end{theorem}

\noindent Clearly,  (\ref{divcond1}) is valid  whenever
(\ref{divcond2}) is satisfied. In light of the comment preceding
Theorem \ref{extradivthm}, we obtain the statement  mentioned in the
abstract as a special case.
\begin{corollary}\label{powdivcor}
Let $\psi$ be any nonnegative arithmetical function. Then $\lambda
(W(\psi))=1$ if there exists an $\epsilon > 0$ such that
\begin{equation}\label{divcond3}
\sum_{n\in\N}\left(\frac{\psi(n)}{n}\right)^{1+\epsilon}\varphi
(n)=\infty \, .
\end{equation}
\end{corollary}

\noindent As pointed out to us by Glyn Harman, Corollary
\ref{powdivcor} is equivalent to Theorem 3.7(iii) in his book
\cite{HarmanMNT}.  
We
will say a little more about the  connection to Harman's work in
\S\ref{thm1pfsec}.
The following Hausdorff dimension  statement is a  consequence of
Corollary \ref{powdivcor} and the recent Mass Transference Principle
\cite{MTP} -- see \S\ref{finalsec} for the corresponding Hausdorff
measure consequence  of Theorem \ref{extradivthm}.
\begin{theorem}\label{thmextradivthm}
Let $\psi$ be any nonnegative arithmetical function. Then  $
\dim(W(\psi))=1 $ if 
\begin{equation} \label{divsumthm2}
\sum_{n\in\N}\left(\frac{\psi(n)}{n}\right)^{1-\epsilon}\varphi
(n) \, = \, \infty \  \qquad \forall \ \epsilon > 0 \,  .
\end{equation}
\end{theorem}

\vspace{1ex}

 It is easily verified that Theorem \ref{thmextradivthm}
implies the following result - see \S\ref{pfcorr1} for the
details.

\begin{corollary}\label{corextradivthm}
Let $\psi$ be any nonnegative arithmetical function. Then  $$
\dim(W(\psi))=1 \qquad {if}   \qquad \sum_{n\in\N} \psi (n)
=\infty \, .
$$
\end{corollary}

\noindent This result has previously been established by R.~C.
Baker and G. Harman -- see Theorem 10.7 in \cite{HarmanMNT}. Our
approach is very different to theirs and enables us to prove the
stronger statement given by Theorem \ref{thmextradivthm} and more
importantly pursue a natural line of questioning that `converges'
to the Duffin-Schaeffer Conjecture -- see \S\ref{finalsec}. To
compare the above dimension statements consider the function $
\psi_{\tau}(n)  :=  n^{-1}{(\log n)^{-\tau}}  $ where $\tau > 1$
is arbitrary.  Regarding Theorem \ref{thmextradivthm}, it is readily
verified that (\ref{divsumthm2}) is satisfied and so $ \dim
(W(\psi_{\tau})) =1 $. On the other hand, $\sum \psi_{\tau}(n) <
\infty $ and so Corollary \ref{corextradivthm} is not applicable.

 The following result is a direct consequence of Corollary \ref{corextradivthm}
and the fact that
$$
\sum_{n\in\N} \frac{\psi (n) }{n} \, \varphi (n) =\infty   \qquad
\Rightarrow \qquad \sum_{n\in\N} \psi (n) =\infty \ .
$$

\begin{corollary}\label{dimDS}
Let $\psi$ be any nonnegative arithmetical function. Then $$
\dim(W(\psi))=1 \qquad {if}   \qquad \sum_{n\in\N}\frac{\psi (n)
}{n}  \, \varphi (n) =\infty \, .
$$
\end{corollary}

\noindent  The corollary can naturally be interpreted as  the
dimension analogue of the Duffin-Schaeffer Conjecture.

\section{Preliminaries for Theorem \ref{extradivthm}}

A well known result due to Gallagher states  that
$\lambda(W(\psi))=0$ or $1$ -- see \cite[Theorem 2.7]{HarmanMNT}.
Furthermore by the Borel-Cantelli Lemma from probability theory
$\lambda(W(\psi))=0$ whenever the sum on the left of
(\ref{divcond1}) is finite. If the sets in the collection
$\{\E_n\}_{n\in\N}$ were pairwise independent then the divergence
part of the same lemma would guarantee that $\lambda (W(\psi))=1$
whenever (\ref{divcond1}) is satisfied. However this turns out not
to be the case and it appears that the best estimate of pairwise
intersection is essentially the following result --  see
\cite[Lemma 2.8]{HarmanMNT} and \cite{VaughanPollington}.
\begin{lemma}
Suppose $\psi$ is a nonnegative arithmetical function and  for
distinct  $m, n\in \N$ let $d:=(m,n)$ and
\[\Delta(m,n):=\max\left\{\frac{\psi (m)}{m},\frac{\psi (n)}{n}\right\}.\]
Then there exists a universal constant $c_1$ with the property
that
\begin{equation}\label{overlap1}
\lambda (\E_m\cap \E_n)\le c_1\lambda (\E_m)\lambda (\E_n)P(m,n),
\end{equation}
where
\[P(m,n):=\prod_{\substack{p|mn/d^2\\p>mn\Delta(m,n)/d}}\left(1-\frac{1}{p}\right)^{-1}.\]
\end{lemma}
The presence of the $P(m,n)$ term cannot be ignored, as it
follows from one of Merten's theorems that
\[\prod_{p|n}\left(1-\frac{1}{p}\right)^{-1}\gg\log\log n\]
for infinitely many $n$. However it is not difficult to show that
there exists a universal constant $c_2$ with the property that
\begin{equation*}
\prod_{\substack{p|n\\p>\sqrt{\log
n}}}\left(1-\frac{1}{p}\right)^{-1}\le c_2
\end{equation*}
for all $n$, so the only time when the $P(m,n)$ term can become
large is when
\begin{equation}\label{pmn1}
\frac{mn\Delta(m,n)}{d}\le \sqrt{\log mn}.
\end{equation}
It is also important to note that for $m\not= n$ if $\E_m$ intersects $\E_n$ then we must have that
\begin{equation*}
\left|\frac{a}{m}-\frac{b}{n}\right|\le \frac{\psi (m)}{m}+\frac{\psi (n)}{n},
\end{equation*}
for some integers $a$ and $b$ with $(a,m)=(b,n)=1$. Multiplying both sides of this equation by $mn$ reveals that
\begin{equation}\label{pmn2}
\lambda (\E_m\cap\E_n)=0\quad\text{ unless }\quad d\le
\psi(m)n+\psi(n)m.
\end{equation}
It stands to reason that conditions (\ref{pmn1}) and (\ref{pmn2})
cannot be satisfied for a large proportion of pairs of integers
$m$ and $n$. So although it is possible for two sets $\E_m$ and
$\E_n$ to have some dependence we wish to show that this cannot
happen on average at the same time that (\ref{divcond2}) is
satisfied. The following well known tool from probability theory
will then suffice to finish the proof -- see \cite[Lemma
2.3]{HarmanMNT}.
\begin{lemma}\label{var1}
Assume that (\ref{divcond1}) is satisfied. Then
\begin{equation*}
\lambda (W(\psi))\ge\limsup_{N\rar\infty}\left(\sum_{n=1}^N\lambda
(\E_n)\right)^2\left(\sum_{m,n=1}^N\lambda (\E_m\cap
\E_n)\right)^{-1}.
\end{equation*}
\end{lemma}

\noindent Finally, we recall the following well known fact:
\begin{equation}
\label{wkf}
 \liminf_{n \to \infty} \frac{\varphi(n) \log \log
n}{n} = e^{-\gamma}  \qquad {\rm where \ } \gamma {\rm \ is \
Euler's \ constant} .
\end{equation}

With all of these observations as a foundation we are ready to
prove Theorem \ref{extradivthm}.

\section{Proof of  Theorem \ref{extradivthm}}\label{thm1pfsec}

Let $\psi$ be a nonnegative arithmetical function with support
$S\subseteq\N$. By appealing to the
Erd\"{o}s-Vaaler Theorem \cite{Va} and to \cite[Theorem
2]{VaughanPollington} we will assume without loss of generality
throughout the proof that $1/n\le \psi (n)\le 1/2$ for all $n\in S$.

\vspace*{2ex}

Our proof is divided into two main steps. In the first step we
work with functions $\psi$ which are essentially constant or zero
on long intervals of the form $[2^{3^K},2^{3^{K+1}})$. This
allows us to take full advantage of the  inequalities (\ref{pmn1}) and
(\ref{pmn2}) stated in the preliminaries.  In the second step we show that we can deal with
general functions $\psi$ by throwing away some of the support and thereby reducing 
the problem back to the `constant' case.  The  extra divergence
condition (\ref{divcond2}) plays a crucial role in both steps.

\vspace*{1ex}

{\em Step 1. \ }For each nonnegative integer $k$ let
\[\psi_k=\min\left\{\psi (n):n\in[2^k,2^{k+1})\cap S\right\}.\] To
begin we will prove the theorem under the extra hypothesis that
\begin{equation}\label{psigrowth1}
\max\left\{\psi (n):n\in\left[2^{3^K},2^{3^{K+1}}\right)\cap
S\right\}\le 2\min\left\{\psi_\ell: 3^K\le \ell < 3^{K+1}\right\},
\end{equation}
for all $K\in\N\cup\{0\}$. Then it is easy to see that we have
\begin{align}\label{E1}
\sum_{n\le 2^{K+1}}\lambda (\E_n)=
\sum_{k=0}^K\sum_{n\in[2^k,2^{k+1})}\frac{\psi (n)\varphi
(n)}{n}\gg\sum_{k=0}^K\frac{\psi_k\cdot|\{n\in S\cap
[2^k,2^{k+1})\}|}{\log\log (2^k)}.
\end{align}

Similarly hypothesis (\ref{divcond2}) guarantees that
\begin{equation}\label{E2}
\sum_{k=0}^\infty2^kf\left(\frac{\psi_k}{2^k}\right)\cdot|\{n\in
S\cap [2^k,2^{k+1})\}|=\infty.
\end{equation}
This clearly ensures that
\begin{equation}\label{setsize1}
|\{n\in S\cap
[2^k,2^{k+1})\}|>\left(k^22^kf\left(\frac{\psi_k}{2^k}\right)\right)^{-1}
\end{equation}
for infinitely many $k$. Furthermore we may throw out all dyadic
blocks for which (\ref{setsize1}) does not hold without affecting
the divergence of (\ref{divcond2}). So without loss of generality we
will now assume that $S$ has been chosen so that for every $k$
either (\ref{setsize1}) holds or else $|S\cap[2^k,2^{k+1})|=0$.

Next notice that the sum
\begin{align*}
\sum_{K=0}^\infty\sum_{n=2^{3^{2K+j}}}^{2^{3^{2K+j+1}}}&f\left(\frac{\psi(n)}{n}\right)\varphi
(n)
\end{align*}
must diverge for $j=0$ or $1$. As it will make no difference to the
proof, we assume that it diverges for $j=0$ and that
\[S\cap \left[2^{3^{2K+1}},2^{3^{2K+2}}\right)=\emptyset\] for all $K$.

Now choose integers $0\le k\le \ell\le 3k$ for which $S_k=S\cap
[2^k,2^{k+1})$ and $S_\ell=S\cap[2^\ell,2^{\ell+1})$ are both
nonempty. For simplicity here and throughout the proof we write
$x=2^k$ and $y=2^\ell.$ If $m\in S_k$ and $n\in S_\ell$ then
\[\frac{\psi_k}{2x}\le\Delta (m,n)\le \frac{2\psi_k}{x},\]
so in order for equations (\ref{pmn1}) and (\ref{pmn2}) to be
satisfied we would have to have that
\begin{equation}\label{gcdcond1}
\frac{y\psi_k}{2 \sqrt{\log 4xy}}\le d\le 2y\psi_k.
\end{equation}
We will establish an upper bound for the number of pairs of integers
$m\in S_k$ and $n\in S_\ell$ which satisfy (\ref{gcdcond1}). First of all we note that
\[\limsup_{n\rar\infty}\frac{\omega (n)\log\log n}{\log n}=1.\]
Thus if we suppose that $0<\epsilon_1 <1/\log 2-1$ is some fixed constant then we may assume that $x_0\in\R$ is chosen large enough so that for all integers $n\ge x_0$ we have
\[\omega (n)\le \frac{(1+\epsilon_1)\log n}{\log\log n}.\]
From this it follows that we can find an $\epsilon_2>0$ such that
\begin{equation}\label{taubound}
\tau (n)\ll \exp\left(\frac{(1-\epsilon_2)\log n}{\log\log n}\right),
\end{equation}
and such that the implied constant is universal. Now suppose that $m\in S_k$ and that $d$ is a divisor of $m$ which satisfies
(\ref{gcdcond1}). Then there at most
\[4\psi_k^{-1}\sqrt{\log 4xy}\]
choices for $n\in S_\ell$ which are divisible by $d$.  Since
$\exp\{(1-\epsilon_2)\log x/\log\log x\}$ grows faster as $x\rar\infty$ than any
power of $\log x$ we have that
\begin{equation} \label{crucial}
\left|\left\{m\in S_k, n\in S_\ell \, : \, \text{(\ref{gcdcond1})
holds}\right\}\right|\ll
|S_k|\psi_k^{-1}\exp\left(\frac{(1-\epsilon_2)\log x}{\log\log
x}\right).
\end{equation}
Now since $\psi_\ell\ge 1/2y$ we use
(\ref{setsize1}) to deduce that
\begin{align}
|S_\ell|&\gg\left(\psi_\ell\cdot\exp\left(\frac{(1/2-\epsilon_2/4)(\log\psi_\ell-\log y)}{\log(\log y-\log\psi_\ell)}\right)\right)^{-1}\nonumber\\
&\gg\psi_k^{-1}\exp\left(\frac{(1-\epsilon_2/2)\log y}{\log\log
y}\right)\label{setsize2}.
\end{align}
Note that the $k^2$ term which appeared in (\ref{setsize1}) has been absorbed by the extra exponent of $\epsilon_2/4$ here. Inequalities (\ref{crucial}) and (\ref{setsize2}) show that
\begin{equation}\label{pairs1}
\left|\left\{m\in S_k, n\in S_\ell \, : \, \text{(\ref{gcdcond1})
holds}\right\}\right|\ll \frac{|S_k||S_\ell|}{(\log\log y)^2}.
\end{equation}

With a view to applying Lemma \ref{var1} we have for each
nonnegative integer $K$ that
\begin{eqnarray}
\sum_{2^{3^{2K}}\le m,n\le 2^{3^{2K+1}}}\lambda (\E_m\cap \E_n) &
\ll & \sum_{3^{2K}\le k\le\ell\le 3^{2K+1}}\sum_{\substack{m\in
S_k\\n\in S_\ell}}\lambda (\E_m)\lambda (\E_n)P(m,n)\nonumber\\
&\ll&\sum_{k,\ell}\left(\sum_{\substack{m\in S_k\\n\in
S_\ell}}'\lambda (\E_m)\lambda (\E_n)+\sum_{\substack{m\in S_k\\n\in
S_\ell}}''\psi(m)\psi(n)\right),\label{var2}
\end{eqnarray}
where $\sum''$ denotes the sum over pairs $m$ and $n$ which satisfy
(\ref{gcdcond1}) and $\sum'$ denotes the sum over pairs which do
not. Also in the summand of $\sum''$ we have used the fact that
\[\frac{\varphi (m)\varphi (n)P(m,n)}{nm}\le 1.\]
By using (\ref{E1}) and (\ref{pairs1}) we find that
\begin{align*}
\sum_{\substack{m\in S_k\\n\in
S_\ell}}''\psi(m)\psi(n)&\ll\psi_k\psi_\ell\cdot\left|\left\{m\in
S_k, n\in S_\ell \, : \, \text{(\ref{gcdcond1})
holds}\right\}\right|\\
&\ll\frac{\psi_k\psi_\ell\cdot|S_k||S_\ell|}{(\log\log x)(\log\log y)}\\
&\ll\left(\sum_{m\in S_k}\lambda (\E_m)\right)\left(\sum_{n\in
S_\ell}\lambda (\E_n)\right).
\end{align*}
Substituting back into (\ref{var2}) this shows that
\begin{equation*}
\sum_{2^{3^{2K}}\le m,n\le 2^{3^{2K+1}}}\lambda (\E_m\cap
\E_n)\ll\left(\sum_{2^{3^{2K}}\le n\le 2^{3^{2K+1}}}\lambda
(\E_n)\right)^2.
\end{equation*}
Finally if $K+1<L$ then
\[\Delta(m,n)\ge (mn)^{-1/2} \text{ whenever } m\in \left[2^{3^{2K}},2^{3^{2K+1}}\right) \text{ and } n\in\left[2^{3^{2L}},2^{3^{2L+1}}\right).\]
In this case the Lemma on page 196 of \cite{VaughanPollington}
ensures that $\lambda(\E_m\cap\E_n)\ll\lambda(\E_m)\lambda(\E_n)$.
Putting this together with our results we have for $K\in\N$ that
\begin{equation*}
\sum_{n,m\le 2^{3^K}}\lambda (\E_m\cap \E_n)\ll\left(\sum_{n\le
2^{3^K}}\lambda (\E_n)\right)^2,
\end{equation*}
and one application of Lemma \ref{var1}, together with Gallagher's
zero-one law, finishes the proof for $\psi$ which satisfy
(\ref{psigrowth1}).

\vspace*{1ex}

{\em Step 2. \ } Let us suppose that $\psi$ is any nonnegative arithmetical
function for which (\ref{divcond2}) holds. We will use $\psi$ to
define a new function $\psi':\N\rar\R$ as follows. For each
$K\in\N\cup\{0\}$, as $n$ runs through the set $S\cap
[2^{3^{2K}},2^{3^{2K+1}})$ the range of values taken by $\psi (n)$
falls in the interval $[2^{-3^{2K+1}},1/2)$. Thus we may choose an
integer $k=k(K)\in\{1,\ldots ,3^{2K+1}-1\}$ so that
\begin{align}
\sum_{\substack{2^{3^{2K}}\le n\le 2^{3^{2K+1}}\\ \psi (n)\in
[2^{-(k+1)},2^{-k})}}f\left(\frac{\psi(n)}{n}\right)\varphi(n)&\ge 3^{-(2K+1)}\sum_{2^{3^{2K}}\le
n\le 2^{3^{2K+1}}}f\left(\frac{\psi(n)}{n}\right)\varphi(n)\nonumber\\
&\gg\sum_{2^{3^{2K}}\le n\le
2^{3^{2K+1}}}f\left(\frac{\psi(n)}{n}\right)\frac{\varphi(n)}{\log n}.\label{divcond4}
\end{align}
Then we set
\[\psi' (n)=\begin{cases}\psi (n) &\text{ if }n\in[2^{3^{2K}},2^{3^{2K+1}})\text{ and }
\psi(n)\in [2^{-(k(K)+1)},2^{-k(K)})\text{ for some }K\in\N,\\0 &\text{ otherwise}.\end{cases}\]
 It is obvious that $W(\psi')\subseteq W(\psi)$. Furthermore, by the ideas used to prove (\ref{setsize2}),
 it follows from (\ref{divcond2}) and (\ref{divcond4}) that
\[\sum_{\substack{n\in\N\\\psi'(n)\not= 0}}\frac{\psi'(n)\varphi(n)}{n}\exp\left(\frac{(1/2-\epsilon_2/8)\log(\psi' (n)/n)}{\log\log (n/\psi'(n))}\right)=\infty ,\] where $\epsilon_2$ is the same as in (\ref{taubound}). Since $\psi'$ satisfies (\ref{psigrowth1}) we can follow almost the exact same argument from the first half of the proof to conclude that $\lambda(W(\psi'))=\lambda(W(\psi))=1$. The only difference is that equation (\ref{setsize1}) will have to be replaced by
\[|\{n\in S\cap
[2^k,2^{k+1})\}|>\left(k^2\psi_k\cdot\exp\left(\frac{(1/2-\epsilon_2/8)\log(\psi_k/x)}{\log\log( x/\psi_k)}\right)\right)^{-1}.\]
However this has no effect on the passage to (\ref{setsize2}) and thus the rest of the proof remains unchanged.

 \hfill $\Box$

\vspace*{4ex}

\noindent {\em Remark.} In the introduction we mentioned that our
Corollary \ref{powdivcor} is equivalent to Theorem 3.7(iii) in
\cite{HarmanMNT}. In fact, the method used to prove Theorem 3.7(iii)
could be used to prove a result which is only slightly weaker than
our Theorem \ref{extradivthm}. Indeed, that this is possible is
explicitly mentioned by Harman in \cite{HarmanAA} in which Theorem
3.7 is first established. However,  with our
 `direct' approach, not only do we get a slightly stronger result but there
 is the added advantage that
the known limitations in Harman's approach -- specifically, see
Theorem 3.2 in \cite{HarmanMNT} and the comments which succeed it --
are no longer necessarily applicable.  In other words,  our approach
can in principle lead to the seriously  stronger and more desirable
statements outlined in \S\ref{finalsec}.   For example,  a sharper upper bound for the quantity
$$\left|\left\{m\in S_k, n\in S_\ell \, : \, \text{(\ref{gcdcond1})
holds}\right\}\right| $$
would lead to stronger results.  The point being highlighted here is  that in our estimate (\ref{crucial}), we do not even make use of the upper bound imposed on $d$ in (\ref{gcdcond1}).

\section{The Mass Transference Principle}

We now turn our attention to establishing Theorem
 \ref{thmextradivthm}. The proof is essentially a nifty application
of the Mass Transference Principle introduced in \cite{MTP}. In a
nutshell, the principle allows us to transfer Lebesgue measure
statements for $\limsup$ sets to Hausdorff measure  statements. The
following version of the Mass Transference Principle has been
simplified and adapted for the particular application we have in
mind. In order to keep the paper self-contained and to introduce
useful and necessary notation, it is convenient to start by defining
Hausdorff measures $\cH^h$.

A {\em dimension function} $h : \R^+ \to \R^+ $ is an increasing,
continuous  function such that $h(r)\to 0$ as $r\to 0 \, $. Suppose
$F$ is  a non--empty subset of $\R$. For $\rho > 0$, a countable
collection $ \left\{B_{i} \right\} $ of balls in $\R$ with radii
$r_i \leq \rho $ for each $i$ such that $F \subset \bigcup_{i} B_{i}
$ is called a {\em $ \rho $-cover}\/ for $F$. Define
$$ \mathcal{
H}^{h}_{\rho} (F) \, := \, \inf  \sum_{i} h(r_i)  \; ,
 $$ where the infimum is over all $\rho$-covers of $F$.  The {\em
Hausdorff $h$--measure}  of $F$ denoted by   $ \mathcal{ H}^{h} (F)$
is defined as
$$
\mathcal{ H}^{h} (F) := \lim_{ \rho \rightarrow 0}
\mathcal{H}^{h}_{\rho} (F)   \; .
$$
In the case that  $h(r) = r^s$ ($s \geq 0$), the measure $ \mathcal{
H}^{h} $ is the more common {\em $s$--dimensional Hausdorff
measure}\/ $\mathcal{ H}^{s} $.  Moreover, the measure
$\mathcal{H}^{1}$  is precisely Lebesgue measure $\lambda$. The
following easy property
$$
\mathcal{ H}^{s}(F)<\infty\quad  \Longrightarrow\quad \mathcal{
H}^{s'}(F)=0\qquad\text{if }s'>s
$$
implies that there is a unique real point $s$ at which the Hausdorff
$s$-measure drops from infinity to zero (unless the set $F$ is
finite so that $\mathcal{ H}^{s}(F)$ is never infinite). This point
is called the {\em Hausdorff dimension} of $F$ and is formally
defined as
$$
\dim F :=  \inf \left\{ s: \mathcal{ H}^{s} (F) =0 \right\} \ = \
\sup \left\{ s: \mathcal{ H}^{s} (F) =\infty \right\} .
$$
Further details regarding Hausdorff measure and dimension can be
found in~\cite{FalcGFS}. We are now in the position to state the
Mass Transference Principle.

Given a dimension function $h$, define the following transformation
on balls in $\R$:
$$
\textstyle B=B(x,r)\mapsto B^h:=B(x,h(r)) \ .
$$
When $h(x)=x^s$ for some $s>0$ we also adopt the notation $B^s$ for
$B^f$. Clearly $B^1=B$. Recall that $\cH^1$ is comparable to the
one-dimensional Lebesgue measure $\lambda$. Given a sequence of
balls $B_i$, $i=1,2,3,\ldots$, as usual its $\limsup$ set  is
$$
\textstyle \limsup_{i\to\infty}B_i:=\bigcap_{j=1}^\infty\
\bigcup_{i\ge j}B_i \ .
$$

\noindent For such limsup sets, the following statement is the key
to obtaining Hausdorff measure statements from Lebesgue statements.

\begin{theorem}[Mass Transference Principle] \label{BVMTP}
Let $\{B_i\}_{i\in\N}$ be a sequence of balls in $\R$ with radii
$r_i \to 0$ as $i\to\infty$.  Let $h$ be a dimension function such
that $x^{-1}h(x)$ is monotonic.  For  any finite  ball $B$ in $\R$,
if
$$
\cH^1\big(\/B\cap\limsup_{i\to\infty}B^h_i{}\,\big)=\cH^1(B) \
$$
then
$$
\cH^h\big(\/B\cap\limsup_{i\to\infty}B^1_i\,\big)=\cH^h(B) \ .
$$
\end{theorem}

For further details and various generalizations of Theorem
\ref{BVMTP} see \cite{MTP, BV06Slicing, BVparis}.

\subsection{Proof of Theorem \ref{thmextradivthm}  \label{pfthm2} }
 Given  $\epsilon > 0 $, consider the function $
\Psi_{\epsilon} $ defined  by
$$
\frac{\Psi_{\epsilon}(n)}{n}  \, := \,
\left(\frac{\psi(n)}{n}\right)^s \qquad {\rm where} \qquad s :=
\frac{1-\epsilon}{1+ \epsilon}  \ .
$$
Then, in view of the divergent sum condition (\ref{divsumthm2}) we
have that
\begin{eqnarray*}
\sum_{n\in\N}\left(\frac{\Psi_{\epsilon}(n)}{n}\right)^{1+\epsilon}\varphi
(n)  \  & = & \
\sum_{n\in\N}\left(\frac{\psi(n)}{n}\right)^{1-\epsilon}\varphi (n)
\ = \ \infty \ .
\end{eqnarray*}

\noindent By Corollary \ref{powdivcor}, it follows that
$$
\cH^1(W(\Psi_{\epsilon})) \, = \, \cH^1(I) \quad {\rm where \ \ }
I:=[0,1) \, .
$$
The set $W(\Psi_{\epsilon})$ is clearly a $\limsup $ set  of balls
centred at reduced rationals $a/n$ with radii $\Psi_{\epsilon}(n)/n
\to 0 $ as $n \to \infty$. On applying the Mass Transference
Principle with $f(r)= r^s$, we immediately deduce that
\begin{equation} \label{ineq}
\cH^s(W(\psi)) \, = \, \cH^s(I) \, = \, \infty \, . \end{equation}
The last equality is valid since $s < 1$.  It follows from the
definition of Hausdorff dimension that
$$
\dim(W(\psi)) \, \geq  \, s \ .
$$
However, $\epsilon >0$ can be made arbitrarily small and so  it
follows that  $ \dim(W(\psi)) \geq  1 $. The complementary upper
bound is a trivial consequence of the fact that $ W(\psi) \subset
\R$.

\hfill $\Box$

\subsection{Proof of Corollary \ref{corextradivthm}  \label{pfcorr1} }

 We are given that $\sum \psi(n) $ diverges. Without loss
of generality, we can assume that
\begin{equation}
\label{wlog}
 \psi(n) \,  < \,  1  \qquad \forall  \ n \in \N \, .
\end{equation}

\noindent Otherwise, we define the function $\psi*$ given by
$\psi^*(n) := \min \{1,\psi(n)\}$. Clearly $  \sum \psi^*(n)$
diverges and since $W(\psi^*) \subset W(\psi) $ the statement
follows on showing that $\dim (W(\psi^*)) =1 $.

\noindent For any  $\epsilon > 0 $, notice that
\begin{eqnarray*}
\sum_{n\in\N}\left(\frac{\psi(n)}{n}\right)^{1-\epsilon}\varphi (n)
\ & \stackrel{(\ref{wlog})}{>}& \ \sum_{n\in\N}\ \psi(n)
\frac{\varphi(n)}{n^{1-\epsilon}}
\\
\ & \stackrel{(\ref{wkf})}{\gg}&  \ \sum_{n\in\N}\ \psi(n)  \ = \
\infty \ .
\end{eqnarray*}

\noindent Thus the divergent sum hypothesis of Theorem
\ref{thmextradivthm} is satisfied and the statement of the corollary
follows.

\hfill $\Box$

\section{Concluding comments \label{finalsec}}

A consequence of the Mass Transference Principle is that the
Duffin-Schaeffer Conjecture (a Lebesgue measure statement) implies
the Generalised Duffin-Schaeffer Conjecture (a Hausdorff measure
statement) -- see \cite{MTP} for the details. In particular, the  Duffin-Schaeffer Conjecture  implies the
following weakening of the generalised conjecture.

\begin{conjecture}
\label{conjA} Let $h$ be a dimension function such that $r^{-1}
h(r)$ is monotonic and suppose that $r^{-1} h(r) \to \infty $ as
$r \to 0 $. Then
$$
\cH^h(W(\psi))= \infty \qquad {if}   \qquad \sum_{n\in\N} \frac{\psi
(n) }{n}  \, \varphi (n) =\infty \, .
$$
\end{conjecture}

\noindent The case when $\cH^h$ is Lebesgue measure (i.e. when
$h(r) =r$) is naturally excluded since $ \cH^h(W(\psi)) \leq
\cH^h([0,1)) = 1 $ and the statement is clearly false. Of course,
in the Lebesgue case the appropriate statement is the
Duffin-Schaeffer Conjecture.

The motivation behind Conjecture \ref{conjA} is to investigate the
Duffin-Schaeffer Conjecture for  measures `arbitrarily' close to
Lebesgue measure. Thus, the divergent sum condition remains
unchanged  but the condition on $h$ means that the measure $\cH^h$
gives `extra weight' to sets compared to Lebesgue measure.
Heuristically, a ball of radius $r$ is given weight $h(r)$ rather
than just $r$. Thus,  Conjecture \ref{conjA} should in principal be
easier to establish than the Duffin-Schaeffer Conjecture. Indeed,
for any $\epsilon
> 0$,   it is easily seen that Corollary \ref{dimDS}  verifies
Conjecture \ref{conjA} for the dimension function $h(r) =
r^{1- \epsilon} $.  Moreover,  let $f$ be as in  Theorem
\ref{extradivthm} and $h$ be any dimension function such that
$r^{-1} h(r)$ is monotonic and
\[\sum_{n\in\N}f\Big(h\Big(\frac{\psi (n)}{n}\Big)\Big)\varphi (n)=\infty.\]
Then Theorem \ref{extradivthm}  together with the Mass Transference
Principle implies  that $\cH^h (W(\psi))=\cH^h([0,1)) $. In
particular, this verifies Conjecture \ref{conjA} for the
dimension function $$ h(r):=r^{1-\frac{1}{1+\log\log 1/r}} \ .
$$ Note that this dimension function tends to zero  as $r$ tends to zero more slowly than $r$
but more quickly than $r^{1-\epsilon}$  for any $\epsilon >0$. Thus,
the next `significant' and natural step towards the conjecture would
be to consider the following  problem.

\begin{problem}
\label{probA} Verify Conjecture \ref{conjA} for the dimension
function $ h(r) := r \log \frac{1}{r} $.
\end{problem}

Related to Conjecture \ref{conjA} is the following weakening of
the Duffin-Schaeffer Conjecture in which the measure $\lambda$
remains unchanged but `extra divergence' is assumed.

\begin{conjecture}
\label{conjB} Let $g$ be an increasing function such that $r^{-1}
f(r) \to 0 $ as $r \to 0 $. Then
$$
 \lambda(W(\psi))= \infty  \qquad {if}   \qquad
\sum_{n\in\N} f \Big( \frac{\psi (n) }{n} \Big)  \, \varphi (n)
=\infty \, .
$$
\end{conjecture}

\noindent In view of the condition imposed on $g$ the above
divergence condition  implies (\ref{divcond1}). Thus,  Conjecture
\ref{conjB} should in principal be easier to establish than the
Duffin-Schaeffer Conjecture.  Indeed, Corollary \ref{powdivcor}
verifies the conjecture for the functions $f(r) = r^{1+ \epsilon} $
where $\epsilon > 0$ is arbitrary. Moreover, Theorem
\ref{extradivthm} verifies Conjecture \ref{conjB} for a function $f$
that  tends to zero as $r$ tends to zero faster than $r \, (- \log
r)^{-\epsilon}$ but more slowly than $r^{1+\epsilon}$ for any
$\epsilon>0$. For this reason we propose that the following problem
represents the next natural step towards Conjecture \ref{conjB}.

\begin{problem}
\label{probB} Verify Conjecture \ref{conjB} for the function $ f(r)
:= r \, \big( \log \frac{1}{r} \big)^{-1} $.
\end{problem}

\noindent On modifying the argument used
to deduce Theorem \ref{thmextradivthm}  from Theorem \ref{extradivthm}
in \S\ref{pfthm2}, it is easily verified that Problem \ref{probB}
implies Problem \ref{probA}. Moreover, and just as easily, we see
that  the Mass Transference Principle enables us to deduce
Conjecture \ref{probA} from Conjecture \ref{probB} whenever $f h (r)
\gg r$ for all $r$ sufficiently small.

\vspace{1ex}

The overall hope is that investigating the weaker conjectures will
yield valuable new insights into the Duffin-Schaeffer Conjecture --
a fundamental  unsolved problems in metric number theory.

\vspace{6ex}

\noindent{\em Acknowledgements. } We would like to thank Glyn Harman
for pointing out the equivalence between our Corollary
\ref{powdivcor} and his Theorem 3.7(iii) in \cite{HarmanMNT}.   SV
would like to thank EPSRC for supporting this research -- in
particular Katharine Bowes who was a great help during the
application stage. Also he would like to thank Fernandez, Hobbs and
Robinson for inspiring the dynamic duo -- Iona and Ayesha -- during
Year 1.

\vspace{10mm}

\noindent Alan K. Haynes: Department of Mathematics, University of
York,

\vspace{0mm}

\noindent\phantom{Alan K. Haynes: }Heslington, York, YO10 5DD,
England.


\noindent\phantom{Alan K. Haynes: }e-mail: akh502@york.ac.uk

\vspace{5mm}

\noindent Andrew D. Pollington: National Science Foundation

\vspace{0mm}

\noindent\phantom{Andrew D. Pollington: }Arlington VA 22230 USA


\noindent\phantom{Andrew D. Pollington: }e-mail: adpolling@nsf.ov

\vspace{5mm}

\noindent Sanju L. Velani: Department of Mathematics, University
of York,

\vspace{0mm}

\noindent\phantom{Sanju L. Velani: }Heslington, York, YO10 5DD,
England.


\noindent\phantom{Sanju L. Velani: }e-mail: slv3@york.ac.uk

\end{document}